\pdfoutput=1
\documentclass{article}
\usepackage{url, latexsym, amssymb, amsmath, amsfonts, amsthm, verbatim, listings, mathrsfs, mathtools, array, multirow, xcolor, xspace, multicol, graphicx, fancyhdr, tcolorbox, setspace, tikz, tikz-cd, indentfirst}

\newcommand{\set}[1]{\left\{ #1 \right\}}
\newcommand{\abs}[1]{\left| #1\right|}
\newcommand{\ar}[1]{\left\langle #1 \right\rangle}
\newcommand{\br}[1]{\left[#1\right]}
\newcommand{\pr}[1]{\left( #1 \right)}

\newcommand{\DS}{\displaystyle}

\newcommand{\eps}{\varepsilon}

\newcommand{\C}{\mathbb{C}}

\newcommand{\g}{\mathfrak{g}}

\newcommand{\h}{\mathfrak{h}}

\renewcommand{\L}{\Lambda}

\newcommand{\Q}{\mathbb{Q}}

\newcommand{\Z}{\mathbb{Z}}
\renewcommand{\phi}{\varphi}
\renewcommand{\>}{\rightarrow}

\newcommand{\ol}{\overline}
\renewcommand{\emph}[1]{\textbf{\textit{#1}}}

\title{Modular Invariance of Characters for Affine Lie Algebras at Subprincipal Admissible Levels}
\author{Victor G. Kac$^1$ and Minoru Wakimoto$^2$}
\usepackage{pdfpages}

\date{%
    $^1$Corresponding author: \texttt{kac@math.mit.edu}\\%
    $^2$Contributing author: \texttt{wakimoto@r6.dion.ne.jp}%
}
\begin{document}
\maketitle

\begin{abstract}
    We prove that the span of normalized characters of subprincipal admissible modules over an affine Lie algebra of subprincipal admissible level $k$ is $SL_2(\Z)$-invariant and find the explicit modular transformation formula.
\end{abstract}

\tableofcontents

\setcounter{section}{-1}

\section{Introduction}

Let $\g$ be an affine Lie algebra over $\C$, and let $L(\Lambda)=L(\Lambda,\g)$ denote the irreducible highest weight module over $\g$ with highest weight $\Lambda$. Recall \cite{KW88} that the weight $\Lambda$ and the corresponding $\g$-module $L(\Lambda)$ are called \textit{admissible} if the following two properties hold:
\begin{align*}
    \ar{\Lambda + \rho, \alpha }\not\in \Z_{\leq 0} \quad \text{for all $\alpha \in \Delta_+^{\mathrm{re}\,\vee}$},\tag{0.1}
\end{align*}
where $\Delta_+^{\mathrm{re}\,\vee}$ is the set of positive real coroots of $\g$;
\begin{align*}
    \Q\set{\alpha \in\Delta_+^{\mathrm{re}\,\vee}\mid \ar{\Lambda + \rho ,\alpha} \in \Z} \supset \Delta_+^{\mathrm{re}\,\vee}.\tag{0.2}
\end{align*}
These $\g$-modules are important since we proved in \cite{KW88} the following formula for their characters $\mathrm{tr}_{L(\Lambda)} e^h,\ h \in \h$:
\begin{align*}
    \prod_{\alpha\in\Delta_+} \pr{1-e^{-\alpha(h)}}^{\mathrm{mult} \ \alpha} \mathrm{tr}_{L(\Lambda)} e^h = \sum_{w\in W^\Lambda} \eps(w) e^{(w.\Lambda)(h)}.\tag{0.3}
\end{align*}
Here $\h$ is the Cartan subalgebra of $\g$, $W^\Lambda$ is the subgroup of the Weyl group $W$ of $\g$, generated by reflection in $\alpha \in \Delta_+^{\mathrm{re}\,\vee}$, \textit{integral} with respect to $\Lambda$, i.e. such that $\ar{\Lambda+\rho,\alpha} \in \Z$,
\begin{align*}
    w. \Lambda = w(\Lambda +\rho) - \rho, \quad w \in W,\tag{0.4}
\end{align*}
is the `shifted' action of $W$, and
\begin{align*}
    \eps(w) = {\det}_\h w.
\end{align*}
Using this formula, we showed that the characters, normalized by a power of $q = e^{-\delta}$ factor, are meromorphic modular functions in the domain $\set{h\in \h\mid \mathrm{Re}\ \delta(h) > 0}$. We also conjectured that these are all irreducible highest weight $\g$-modules $L(\Lambda)$ with this property.

From the viewpoint of the theory of vertex algebras the most interesting case is that of a non-twisted $\g$, associated to a simple Lie algebra $\ol{\g}$, and $\Lambda = k \Lambda_0$, where $k \in \C$, when $L(k \Lambda_0)$ admits a structure of simple affine vertex algebra $V_k(\ol{\g})$ of level $k$. If $k \Lambda_0$ is an admissible weight, $k$ is called an \textit{admissible level}. It is easy to deduce from the classification of admissible weights in \cite{KW89} that $k$ is an admissible level only when it is a rational number, and there are only the following two possibilities \cite{KW08}: 

(principal $k$) $k = -h^\vee + \frac{p}{u}$, where $p$ and $u$ are coprime positive integers, $u$ is coprime with the lacety $r^\vee$ of $\ol{\g}$, and $p \geq h^\vee$, the dual Coxeter number of $\g$; 

(subprincipal $k$) $k=-h^\vee + \frac{p}{u}$, where $p$ and $u$ are coprime positve integers, $u$ is divisible by $r^\vee > 1$, and $p\geq h$, the Coxeter number of $\g$.

The above mentioned conjecture leads to 
\vspace{0.3cm}\\
\noindent \textbf{Conjecture 0.1.} The series $\mathrm{tr}_{V_k(\g)} e^{2\pi i \pr{L_0 - \frac{C_k}{24}}\tau}$ converges to a holomorphic modular function in the upper half-plane Im $\tau > 0$ if and only if $k$ is either principal or subprincipal admissible level. Here $C_k = \frac{k \dim \ol{\g}}{k+h^\vee}$.
\vspace{0.5cm}\\
\indent Note that the `if' part holds due to \cite[formula (52)]{KW25} if $k$ is a principal admissible level, and due to the following similar formula for $k=-h^\vee + \frac{p}{u}$ subprincipal admissible level: 
\begin{align*}
    \mathrm{tr}_{L(\Lambda)} q ^{L_0 - \frac{C_k}{24}} = \frac{\sum_{\gamma \in \ol{Q}} d(\ol{\Lambda}+p\gamma) q^{\frac{u}{2p}\abs{\ol{\Lambda}+\ol{\rho} + p\gamma}^2}}{\eta(\tau)^{\dim \ol{\g}}},\tag{0.5}
\end{align*}
if $\Lambda$ is a subprincipal admissible weight of level $k$ whose restriction to $\ol{\h}$ is dominant integral for $\ol{\g}$ where
\begin{align*}
    d(\beta) = \prod_{\alpha\in \ol{\Delta}_+} \frac{\pr{\beta + \ol{\rho} \mid \alpha}}{\pr{\ol{\rho} \mid \alpha}}, \quad \beta \in \ol{\h}^*.
\end{align*}
In particular, (0.5) holds for $\Lambda = k \Lambda_0$.

Note that (0.5) is a modular function since it is a specialization of the normalized character of $L(\Lambda)$.

The most important among the admissible weights $\Lambda$ of principal admissible level $k$ is the set $P_{+,\mathrm{pr}}^k$ of \textit{principal admissible weights}. Such $\Lambda$ are defined by the property that the set of real coroots, integral with respect to $\Lambda$, is isometric to that with respect to $k\Lambda_0$, or, equivalently, its set of simple roots is isometric to the set 
\begin{align*}
    S_{(u)}^\vee = \set{u K - \theta^\vee, \alpha_1^\vee, \cdots ,\alpha_l^\vee},\tag{0.6}
\end{align*}
where $K$ is the canonical central element of $\g$, $\theta$ is the highest root of $\ol{\g}$, and $\alpha_1^\vee, \cdots ,\alpha_l^\vee$ are simple coroots of $\ol{\g}$. (If $\g$ is of type $A_\ell^{(1)}$, then $P_{+,\mathrm{pr}}^k$ is the set of all admissible weights, but not for other types.)

It is shown in \cite{KW89} that $P_{+,\mathrm{pr}}^k$ consists of weights of the form
\begin{align*}
    \Lambda = \pr{t_\beta \ol{y}} . \pr{\Lambda^0 + \frac{1-u}{u} p \Lambda_0},\tag{0.7}
\end{align*}
where $\beta \in \ol{P}$ and $y\in\ol{W}$ are such that
\begin{align*}
    t_\beta \ol{y} \pr{S_{(u)}^\vee} \subset \Delta_+^{\mathrm{re}\,\vee},
\end{align*}
and $\Lambda^0$ is a dominant integral weight for $\g$ of level $p-h^\vee$.

The importance of principal admissible weights is due to the fact, proved in \cite{KW89}, that the span of normalized characters of $\g$-modules $L(\Lambda)$ with $\Lambda \in P_{+,\mathrm{pr}}^k$ is $SL_2(\Z)$-invariant, and the fact, proved in \cite{A16} that all $\g$-modules $L(\Lambda)$, which descend to $V_k(\g)$ are of the form (0.7), with $\beta = 0, \ol{y} =e$, i.e. $\Lambda$ is quasidominant.

The set of admissible weights of subprincipal admissible level $k$, denoted by $P_+^k$ in the paper, is defined in a similar way. A weight of such level $k$ is called \textit{subprincipal admissible} if the set of coroots of $\g$, which are integral with respect to $\Lambda$, is isometric to that for $k\Lambda_0$, i.e. to the set
\begin{align*}
    S_{(u)}^{\#\vee} = \set{uK - \theta_s^\vee, \alpha_1^\vee, \cdots ,\alpha_l^\vee},\tag{0.8}
\end{align*}
where $\theta_s$ is the highest short root of $\ol{\g}$ (cf. (0.6)).

In this case the set of coroots of $\g$, integral with respect to $\Lambda$, is isometric to the set of coroots of the twisted affine Lie algebra $\g^\# = \g\pr{A^\#}$, whose Cartan matrix $A^\#$ is $A^{\top'}$, where $A$ is the Cartan matrix of $\g$, $\top$ stands for transpose and $'$ stands for taking the adjacent Cartan matrix \cite[\S 13.9]{K90}.
Recall that $A'=A$ for all affine Cartan matrices, except that $'$ switches $A_{2\ell -1}^{(2)}$ with $D_{\ell + 1}^{(2)}$. 

The description of the set of subprincipal admissible weights is similar to that of principal ones, see Theorem 2.1.
There is an important difference, however. All $k\in \Z_{\geq 0}$ are principal admissible levels and the corresponding principal admissible weights are just level $k$ dominant integral weights, hence the span of the normalized characters of $\g$-modules $L(\Lambda)$ with principal admissible $\Lambda$ of level $k\in \Z_{\geq 0}$ is $SL_2(\Z)$-invariant \cite{K90}. This is not the case for subprincipal admissible levels, which cannot be integer numbers. It is related to the fact that the numerators of characters of subprincipal admissible $L(\Lambda)$ are expressed in terms of the Weyl group $W^\#$ of the twisted affine Lie algebra $\g^\#$. But the span of normalized characters of modules $L(\Lambda, \g^\#)$ with $\Lambda$ dominant integral is invariant only with respect to the action of the subgroup $\Gamma_0(r^\vee)$ of $SL_2(\Z)$ \cite[Theorem 13.9]{K90}. 

It is quite surprising that nevertheless the span of the normalized characters of subprincipal admissible $\g$-modules $L(\Lambda)$, $\Lambda\in P_+^k$, is $SL_2(\Z)$-invariant. This is our main Theorem 5.1. Due to the character formula (3.6), it suffices to prove $SL_2(\Z)$-invariance of the span of the numerators $A_{\Lambda + \rho}$, since $\C R_\g$, where $R_\g$ is the Weyl denominator, is $SL_2(\Z)$-invariant. The $SL_2(\Z)$-invariance of the span of $\set{A_{\Lambda + \rho}}_{\Lambda \in P_+^k}$ is proved, using their expression in terms of Jacobi forms in Proposition 5.1.

Note that Theorem 5.1 implies the same modular transformation formulas for quantum Hamiltonian reduction $H_f(\Lambda)$ for subprincipal admissible $\Lambda$, as those, given by \cite[Theorem 4]{KW25} for the principal admissible ones.

In the last section of the paper we show that the simple vertex algebra $V_k(\ol{\g})$ of subprincipal admissible level $k$ is modular invariant in the sense of \cite{KW25}, where this is proved for principal admissible levels. We also find the conditions, for which the quantum Hamiltonian reduction $H_f$ at subprincipal admissible level $k$ is a modular invariant vertex algebra.

For the basics of the theory of affine Lie algebras and their representations, and of the theory of Jacobi forms, see the book \cite{K90}.

\section{Recollection on affine Lie algebras and their root systems}

Let $\ol{\g}$ be a simple finite-dimensional Lie algebra over $\C$ of rank $\ell$ and lacety $r^\vee > 1$, i.e. $\g$ is of type $B_\ell, C_\ell, F_4,$ or $G_2$, and $r^\vee=2,2,2,$ or 3 respectively. Choose a Cartan subalgebra $\ol{\h}$ of $\ol{\g}$. Choose a set of simple roots $\ol{\Pi} = \set{\alpha_1,\cdots,\alpha_\ell}$, normalize the invariant bilinear form on $\ol{\g}$ by the condition $(\theta\mid \theta) = 2$, where $\theta$ is the highest root. We identify $\ol{\h}$ with $\ol{\h}^*$, using this bilinear form. We denote by $\theta_s$ the highest short root. Then
\begin{align*}
    \pr{\theta_s \mid \theta_s} = \frac{2}{r^\vee}.\tag{1.1}
\end{align*}
Let $\ol{\Delta}$ be the set of all roots of $\ol{\g}$, then $\ol{\Delta} = \ol{\Delta}_\ell \cup \ol{\Delta}_s$ is a union of the sets of long and short roots $\alpha$, so that $(\alpha\mid \alpha) = 2$ and $\frac{2}{r^\vee}$ respectively. As usual, for a root $\alpha \in \ol{\Delta}$ we denote by $\alpha^\vee = \frac{2\alpha}{(\alpha\mid \alpha)}$ the corresponding coroot.

Let $\ol{W}$ be the Weyl group of $\ol{\g}$, let $\ol{Q} = \sum_{i=1}^\ell \Z \alpha_i, \ \ol{Q}^\vee = \sum_{i=1}^\ell \Z \alpha_i^\vee$ be the root and coroot lattices, and let $\ol{P}^\vee$ and $\ol{P}$ be their dual lattices respectively. We have
\begin{align*}
    r^\vee \ol{Q} \subset \ol{P}^\vee.\tag{1.2}
\end{align*}
Let $\ol{\rho} \in \ol{\h}^*$ be the Weyl vector for $\ol{\g}$, i.e. $\pr{\ol{\rho} \mid \alpha_i^\vee} = 1$ for $i=1,\cdots,\ell$. Then
\begin{align*}
    r^\vee \ol{\rho} \in \ol{P}^\vee.\tag{1.3}
\end{align*}

Let $\g$ be the affine Lie algebra, associated to $\ol{\g}$ (see \cite{K90} for details). Recall that it is a Kac-Moody algebra $\g(A)$ with Cartan matrix $A$ of type $B_\ell^{(1)}, C_\ell^{(1)}, F_4^{(1)},$ or $G_2^{(1)}$ respectively. Let $\g\pr{A^\#}$ be the twisted affine Lie algebra with the Cartan matrix $A^\#$ of type $D_{\ell+1}^{(2)}, A_{2\ell -1}^{(2)}, E_6^{(2)}, $ or $D_4^{(3)}$ respectively.

Recall that $\g$ can be explicitly construced as $\g = \ol{\g}[t,t^{-1}] \oplus \C K \oplus \C d$, where $K$ is a central element, and the remaining brackets are $(a,b\in \g, \ m,n\in Z)$:
\begin{align*}
    \br{at^m, bt^n} = [a,b]t^{m+n} + m \delta_{m,-n}(a\mid b) K, \ [d,t^ma]=mt^m a.\tag{1.4}
\end{align*}
The invariant symmetric non-degenerate bilinear form on $\ol{\g}$ extends to that on $\g$ by letting
\begin{gather*}
    \pr{at^m \mid bt^n} = (a\mid b)\delta_{m,-n},\quad (at^m \mid \C K + \C d) = 0,\\ (K\mid d) = 1, \quad (K\mid K) = (d\mid d) = 0.\tag{1.5}
\end{gather*}

Let $\h=\C d + \ol{\h} + \C K$ be the Cartan subalgebra of $\g$. We identify it with $\h^*$ using the bilinear form $(\bullet \mid \bullet)$. Then $\delta \in \h^*$, defined by
\begin{align*}
    \delta \mid_{\ol{\h} + \C K} = 0,\quad \delta(d) = 1,\tag{1.6}
\end{align*}
is identified with $K$. Let $\Lambda_0\in\h^*$ be the element, identified with $d$. The set of roots of $\g$ with respect to $\h$ is $\Delta = \Delta^\mathrm{re} \cup \Delta^\mathrm{im}$, where
\begin{align*}
    \Delta^\mathrm{re} = \set{\alpha + n\delta \mid \alpha \in \ol{\Delta}, n \in \Z}, \quad \Delta^\mathrm{im} = \set{n\delta,\, n \in \Z\setminus\set{0}}.\tag{1.7}
\end{align*}
The set of simple roots is 
\begin{align*}
     \Pi = \set{\alpha_0 = \delta- \theta} \cup \ol{\Pi}.\tag{1.8}
\end{align*}
The multiplicities of roots in $\Delta^\mathrm{re}$ (resp. $\Delta^\mathrm{im}$) are 1 (resp. $\ell$).

Recall the translation operator $t_\gamma$ on $\h$ for $\gamma \in \ol{\h}^*$:
\begin{align*}
    t_\gamma(h) = h + \delta(h) \gamma - \pr{\frac{(\gamma \mid \gamma)}{2}\delta(h) + \gamma(h)}K,\quad h\in \h.\tag{1.9}
\end{align*}
For a subgroup $S\subset \ol{\h}^*$ let $t_S =\set{t_\gamma \mid \gamma \in S}$, which is an abelian group. Then the Weyl group $W$ of $\g$ is a semidirect product of $\ol{W}$ and the abelian normal subgroup $t_{\ol{Q}^\vee}$.

We coordinatize $\h$ by writing $h\in\h$ in the form $h=2\pi i(-\tau d + z + tK)$, where $\tau, t \in \C$, $z\in \ol{\h}$, and we write a function $f$ on $\h$ in the form $f = f(\tau, z, t)$; we also write $f(\tau, z) = f(\tau, z, 0)$.

The main character of this paper is the twisted affine Lie algebra $\g^\# = \g\pr{A^\#}$. We identify its Cartan subalgebra with the Cartan subalgebra $\h$ of $\g$. The set of roots of $\g^\#$ is $\Delta^\#= \Delta^{\#\mathrm{re}}\cup \Delta^{\#\mathrm{im}}$, where
\begin{gather*}
    \Delta^{\#\mathrm{re}} = \set{\alpha + nr^\vee \delta \mid \alpha \in \ol{\Delta}_\ell , n \in \Z} \cup \set{\alpha + n\delta \mid \alpha \in \ol{\Delta}_s, n \in \Z},\\
    \Delta^{\#\mathrm{im}} = \set{n\delta \mid n \in \Z \setminus\set{0}}.\tag{1.10}
\end{gather*}
Multiplicities of roots from $\Delta^{\#\mathrm{re}}$ are 1, but of the roots from $\Delta^{\#\mathrm{im}}$ are different:
\begin{align*}
    \mathrm{mult}\ n\delta = \frac{N-\ell}{r^\vee -1} \text{ (resp. $\ell$) if } n\in r^\vee \Z \text{ (resp. $n\not\in r^\vee \Z$)},\tag{1.11}
\end{align*}
where $A^\#$ is of type $X_N^{(r^\vee)}$.

The set of simple roots for $\g^\#$ is 
\begin{align*}
    \Pi^\# = \set{\delta - \theta_s} \cup \ol{\Pi}.\tag{1.12}
\end{align*}
The Weyl group $W^\#$ of $\g^\#$ is
\begin{align*}
    W^\# = \ol{W} \ltimes t_{r^\vee\ol{Q}}.\tag{1.13}
\end{align*}
As before, for a non-isotropic root $\alpha$, we let $\alpha^\vee = \frac{2\alpha}{(\alpha \mid \alpha)}$ be the corresponding coroot. Then
\begin{align*}
    \Delta^{\#\mathrm{re}\vee} = \set{nK^\# + \alpha^\vee \mid \alpha^\vee \in \ol{\Delta}^\vee, n \in \Z},\tag{1.14}
\end{align*}
and its subset of simple coroots is
\begin{align*}
    \Pi^{\#\vee}= \set{K^\# - \theta_s^\vee , \alpha_1^\vee, \cdots, \alpha_\ell^\vee},\tag{1.15}
\end{align*}
where
\begin{align*}
    K^\# = r^\vee K.\tag{1.16}
\end{align*}

Let $\rho^\#$ be the Weyl vector for $\g^\#$, i.e.
\begin{align*}
    \pr{\rho^\# \mid \alpha} = 1 \qquad \text{for all $\alpha \in \Pi^{\# \vee}$.}\tag{1.17}
\end{align*}

\section{Subprincipal admissible weights}

Recall that a \textit{subprincipal admissible level} for the affine Lie algebra $\g$ is a rational number $k$ of the form
\begin{align*}
    k = \frac{p}{u} - h^\vee,\quad\text{where $p\in \Z_{\geq h}, u = u' r^\vee , u' \in \Z_{>0}, \gcd(p,u)=1,$}\tag{2.1}
\end{align*}
where $h$ is the Coxeter number of $\g$. Recall that
\begin{align*}
    h = h^{\#\vee}, \quad\text{the dual Coxeter number of $\g^\#$.}\tag{2.2}
\end{align*}
Let
\begin{align*}
    S_{(u)}^{\# \vee} = \set{u K - \theta_s^\vee, \alpha_1^\vee, \cdots ,\alpha_\ell^\vee}.\tag{2.3}
\end{align*}
The subgroup $W^{\#}_{(u)}$ of $W^\#$, generated by reflection in these coroots, is
\begin{align*}
    W_{(u)}^{\#} = \ol{W} \ltimes t_{u\ol{Q}}.\tag{2.4}
\end{align*}

\noindent \textbf{Lemma 2.1.}
\begin{enumerate}
    \item[(a)] For each $\beta \in \ol{P}^\vee$ there exist unique $\ol{y}\in \ol{W}$ and $\gamma \in \ol{Q}$, such that
    \begin{align*}
        {t_{\beta+ u \gamma} }\ol{y} S_{(u)}^{\# \vee} \subset \Delta_+^{\mathrm{re}\vee}.\tag{2.6}
    \end{align*}
    \item[(b)] Let $\Phi_{(u)} = \set{y\in\ol{W}\ltimes t_{\ol{P}^\vee} \mid y\pr{S_{(u)}^{\#\vee}}\subset \Delta_+^{\mathrm{re}\vee}}.$ Then (2.6) induces a bijective map
    \begin{align*}
        \phi: \ol{P}^\vee/u\ol{Q} \xrightarrow[]{\sim} \Phi_{(u)}\tag{2.7}
    \end{align*}
\end{enumerate}
\textit{Proof.} Let $\gamma_0 = u {K} - \theta_s^\vee$, $\gamma_i = \alpha_i^\vee$ for $i=1,\cdots,\ell$, to simplify notation. In order to prove (a), if suffices to show that, given $\beta\in\ol{P}^\vee$ there exist a unique $\ol{y}\in \ol{W}$ and a unique $\gamma\in \ol{Q}$, satisfying
\begin{align*}
    \pr{t_{\beta+u\gamma} \ol{y} (\gamma_i) \mid \frac{1}{h}\rho^\vee}>0\quad\text{for $i=0,\cdots,\ell.$}\tag{2.8}
\end{align*}
Indeed, this would imply that the element $\ol{y}^{-1}t_{-u\gamma}t_{-\beta}\pr{\frac{1}{h}\rho^\vee}$ lies in the fundamental Weyl chamber with respect to $S_{(u)}^{\#\vee}$.

Note that
\begin{align*}
    t_{-\beta}\pr{\frac{1}{h}\rho^\vee} \equiv \frac{1}{h}\rho^\vee - \beta\text{ mod }\C\delta,\tag{2.9}
\end{align*}
Since $(\beta \mid \alpha)\in \Z$ and $\pr{\frac{1}{h}\rho^\vee\mid \alpha}\not \in \Z$ for all $\alpha\in\ol{\Delta}$, we see from (2.9) that
\begin{align*}
    \pr{\frac{1}{h}\rho^\vee\mid n\delta + \alpha} \not\in \Z\quad\text{for all $\alpha\in\ol{\Delta},n\in \Z$.}
\end{align*}
Since $\pr{\beta \mid n \delta + \alpha} \in \Z$ for all $\alpha\in\ol{\Delta}, n\in\Z$, we conclude that
\begin{align*}
    \pr{\frac{1}{h}\rho^\vee-\beta\mid \alpha}\not\in \Z \quad\text{for all $\alpha\in\Delta^\mathrm{re}$.}
\end{align*}

This implies that $t_{-\beta}\pr{\frac1h\rho^\vee}$ is a regular element with respect to $S_{(u)}^{\#\vee}$, hence there exists a unique $w\in W_{(u)}^{\#}$, such that $w\pr{t_{-\beta}\pr{\frac1h \rho^\vee}}$ lies in the fundamental Weyl chamber for $S_{(u)}^{\#\vee}$. Writing this element as $w=\ol{y}^{-1} t_{-u\gamma}$, we obtain (2.8).

In order to prove (b), note that we just constructed a well-defined map
\begin{align*}
    \phi:\ol{P}^\vee/u\ol{Q} \> \Phi_{(u)},\quad \beta \mapsto t_{\beta+u\gamma} \ol{y}.
\end{align*}
The map $\phi$ is surjective since if $y=t_\beta \ol{y} \in t_{\ol{P}^\vee} \rtimes \ol{W}$, such that $t_\beta \ol{y}\pr{S_{(u)}^{\#\vee}}\subset \Delta_+^\vee$, then $\phi(\beta) = y$ by definition of $\phi$.

To prove injectivity of $\phi$, suppose that we have $\beta,\beta' \in \ol{P}^\vee$, $\gamma,\gamma' \in \ol{Q}$ and $\ol{y},\ol{y}' \in \ol{W}$, such that $t_{\beta + u\gamma} \ol{y} = t_{\beta' + u\gamma'}\ol{y}'$. Then
\begin{align*}
    t_{\beta + u \gamma - (\beta' + u \gamma')} = \ol{y}' \ol{y}^{-1} \in t_{\ol{P}^\vee} \cap \ol{W} = \set{e},
\end{align*}
hence $\beta + u\gamma = \beta' + u\gamma'$, $(\beta-\beta') = u(\gamma'-\gamma) \in u\ol{Q}$, i.e. $\beta\equiv \beta' \text{ mod } u\ol{Q}$.\qed
\vspace{0.4cm}\\
\noindent \textbf{Remark 2.1.} The analogue of Lemma 2.1 in the principal admissible case (proved in the same way) is as follows: For $u\in \Z_{\geq 1}$, such that $\gcd\pr{u,r^\vee}=1$, let
\begin{align*}
    S_{(u)}^\vee = \set{uK - \theta^\vee, \alpha_1^\vee, \cdots, \alpha_\ell^\vee}.
\end{align*}
Then
\begin{enumerate}
    \item[(a)] For $\beta \in \ol{P}^\vee$, there exist unique $\ol{y} \in \ol{W}$ and $\gamma \in \ol{Q}^\vee$, such that $t_{\beta + u\gamma} \ol{y} \pr{S_{(u)}^\vee}\subset \Delta_+^{\mathrm{re}\vee}$.
    \item[(b)] The map $\phi^\mathrm{pr} : \ol{P}^\vee / u\ol{Q}^\vee \> \set{y\in\ol{W} \ltimes t_{\ol{P}^\vee} \mid y\pr{S_{(u)}^\vee}\subset \Delta_+^\vee}$ is well defined and bijective.
\end{enumerate}

Let $k$ be a subprincipal admissible level as in (2.1), with the numerator $p$ and the denominator $u=r^\vee u'$. An admissible weight $\Lambda$ is called \textit{subprincipal admissible} for $\g$ if the set of real coroots $\set{\alpha\in\Delta_+^{\vee \mathrm{re}}\mid (\Lambda \mid \alpha)\in \Z}$ has the set of simple roots isometric to $S_{(u)}^{\#\vee}$. Denote by $P_+^k$ the set of subprincipal admissible weights of level $k$.

Let
\begin{align*}
    P_{+,e}^k &= \left\{\Lambda \in \h^* \bigg{|}\ \Lambda + \rho \in \frac{p}{u}\Lambda_0 + \ol{P}, (\Lambda+\rho\mid u \delta - \theta_s^\vee)>0,\right.\\
    &\quad \left.\phantom{\bigg{|}}(\Lambda +\rho\mid\alpha_i^\vee) > 0 \text{ for $i=1,\cdots,\ell$}\right\},\tag{2.10}\\
    P_{+,y}^k &= \set{y.\Lambda \mid  \Lambda \in P_{+,e}^k}\text{ for $ y \in \Phi_{(u)}$}.\tag{2.11}
\end{align*}
In the same way as for the principal admissible weights in \cite[Theorem 2.1]{KW89} one can prove the following theorem.
\vspace{0.4cm}\\
\noindent \textbf{Theorem 2.1.} $\DS P_+^k = \bigcup_{y\in \Phi_{(u)}} P_{+,y}^k$, where
\begin{align*}
    P_{+,y}^k = \set{y\pr{\Lambda^0 + \rho^\# + \frac{p(1-u')}{u}\Lambda_0} - \rho\ \bigg{|}\ \Lambda^0 \in P_+^{\#(p-h)}},
\end{align*}
$P_+^{\# (p -h)}$ denotes the set of all dominant integral weights of level $p-h$ for $\g^\#$, and $h$ is the Coxeter number of $\g$.

\section{Normalized characters of subprincipal admissible highest weight $\g$-modules}

For $\Lambda \in \h^*$ of level $k\neq - h^\vee$ let $L(\Lambda, \g)$ denote the irreducible highest weight $\g$-module with highest weight $\Lambda$. Let 
\begin{align*}
    m_\Lambda = \frac{|\Lambda + \rho|^2}{2(k+h^\vee)}, \quad \text{where $\rho =\ol{\rho} + h^\vee \Lambda_0$,}\tag{3.1}
\end{align*}
and let
\begin{align*}
    R_\g = q^{m_0} \prod_{\alpha \in \Delta_+} \pr{1- e^{-\alpha}}^{\mathrm{mult} \ \alpha}\tag{3.2}
\end{align*}
be the normalized Weyl denominator for $\g$. Hereafter
\begin{align*}
    q = e^{-\delta}.\tag{3.3}
\end{align*}
Recall that the normalized character of $L(\Lambda)$ is defined as
\begin{align*}
    \mathrm{ch} L(\L,\g) = q^{m_\L - m_0}\mathrm{tr}_{L(\L,\g)} e^{h}, \quad h \in \h^*.\tag{3.4}
\end{align*}
For a subprincipal admissible weight $\L\in P_{+,y}^k$, let
\begin{align*}
    A_{\L + \rho} = q^{m_\L}\sum_{w\in y W^\#_{(u)}y^{-1}}\eps (w) e^{w(\L + \rho)}.
\end{align*}
Theorem 3.1 from \cite{KW89} gives the following formula
\begin{align*}
    A_{\L + \rho} = R_\g \mathrm{ch} L(\L,\g) .\tag{3.6}
\end{align*}
Let $\g_{(u),y}^\#$ be the affine subalgebra of $\g^\#$, with the set of simple coroots $y\pr{S^{\# \vee}_{(u)}}$, and let $R_{(u),y}^\#$ be its Weyl denominator. Let $\rho_u^\#$ be the Weyl vector for $\g_{(u),e}^\#$; it is given by 
\begin{align*}
    \rho_u^\# = \frac{h}{u'}\L_0+\ol{\rho}.\tag{3.7}
\end{align*}
We write
\begin{align*}
    \L + \rho = \L^{\# 0} + \rho_u ^\#.\tag{3.8}
\end{align*}

\noindent\textbf{Lemma 3.2.} Let $\L \in P_{+,y}^k$ be a subprincipal admissible weight, where $y \in \Phi_{(u)}$, so that $\L = y(\L_e + \rho) - \rho$, where $\L_e \in P_{+,e}^k$, and let $\L_y^{\# 0} = y\pr{(\L_e)^{\# 0}}$. Then
\begin{itemize}
    \item[(a)] $\L_y^{\# 0}$ is a dominant integral weight for $\g_{(u),y}^{\#}$.
    \item[(b)] $A_{\L + \rho} = R_{(u), y}^\# \mathrm{ch}\pr{L \pr{\L_y^{\# 0}, \g_{(u),y}^{\#} }}$.
\end{itemize}

\noindent \textit{Proof.} It suffices to prove (a) for $y=e$. In this case
\begin{align*}
    \L + \rho = \frac{p}{u}\L_0 + \ol{\L} + \ol{\rho},\text{ where $\ol{\L} \in \ol{P}$.}\tag{3.9}
\end{align*}
Since $\rho $ satisfies
\begin{align*}
    \rho(K) = h^\vee, \ \rho(\alpha_i^\vee) = 1 \text{ for $i=1,\cdots, \ell$,}\tag{3.10}
\end{align*}
we obtain
\begin{align*}
    (\L +\rho) (uK - \theta_s^\vee) = u \cdot \frac{p}{u} - (\L +\rho)(\theta_s^\vee)\in \Z.\tag{3.11}
\end{align*}
Formula (3.11) together with (2.11) gives
\begin{align*}
    (\L +\rho) (uK - \theta_s^\vee) \in \Z_{>0}, \ (\L +\rho) (\alpha_i^\vee )\in\Z_{>0},\text{ for $i=1,\cdots,\ell$.}\tag{3.12}
\end{align*}
Hence
\begin{align*}
    (\L +\rho) (\gamma) \in \Z_{>0}\text{ for all $\gamma \in S_{(u)}^{\#\vee}$.}\tag{3.13}
\end{align*}
But $\L +\rho = \L^{\# 0} + \rho_u^\#$, where $\rho_u^\#$ is the Weyl vector for $\g_u^\#$, so that $\rho_u^\#(\gamma) = 1$ for all $\gamma \in S_{(u)}^{\#\vee}$.

Hence, by (3.13), we obtain
\begin{align*}
    \L^{\#0}(\gamma) \in \Z_{\geq 0} \text{ for all $\gamma \in S_{(u)} ^{\#\vee}$},
\end{align*}
proving (a).

(b) Since $\L_y^{\# 0}$ is a dominant integral weight for $\g_{(u),y}^\#$ and $y W_{(u)}^{\#} y^{-1}$ is the Weyl group for $\g_{(u),y}^\#$, by Theorem 3.1 of \cite{KW89}, we have
\begin{align*}
    R_{(u),y}^\# \mathrm{ch}\  L(\L_y^{\#0},\g_{(u),y}^\#)  = \pr{q^{u'}}^{\frac{|\L_y^{\#0}+y \rho_u^\#|^2}{2\pr{k_u^\# + h^{\# \vee}}}}\sum_{w \in y W_{(u)}^{\#} y^{-1}} \eps(w) e^{w\pr{y(  \L_e^{\#0}+ \rho_u^\#)}},\tag{3.14}
\end{align*}
where $R_{(u),y}^\#$ is the Weyl denominator for $\g_{(u),y}^\#$, and
\begin{align*}
    k_u^\# + h^{\# \vee} = \pr{ y(\L_e^{\# 0} + \rho_u^\#)\mid u'K} = u'(k + h^\vee)\tag{3.15}
\end{align*}
since $y(\L_e^{\# 0} + \rho_u^\#) = \L + \rho$.

Hence the normalizing factor in the RHS of (3.14) equals $q^{m_\L}$, and the sum factor equals to $\sum_{w \in y W_{(u)}^{\#} y^{-1}} \eps(w) e^{w(\L + \rho)}$. Hence (b) follows from (3.14).\qed

\section{Theta forms}
Let $k=\frac{p}{u} - h^\vee$ be a subprincipal admissible level with the denominator $u=u'r^\vee$, see (2.1). Let
\begin{align*}
    n = \frac{p}{r^\vee} \quad \pr{\in \frac{1}{r^\vee}\Z_{>0}}.\tag{4.1}
\end{align*}
Consider the weights of the form
\begin{align*}
    \lambda = \frac{n}{u'}\L_0 + \ol{\lambda}, \quad \text{where $\ol{\lambda}\in \pr{r^\vee \ol{Q}}^* = \frac{1}{r^\vee}\ol{P}^\vee$}.\tag{4.2}
\end{align*}
In order to study modular invariance of subprincipal admissible characters, introduce the following Jacobi theta forms:
\begin{align*}
    \Theta_\lambda^{\#(u)} = q^{\frac{(\lambda \mid \lambda) u'}{2n}}\sum_{\gamma \in r^\vee \ol{Q}} e^{t_{u'\gamma}(\lambda)} = \sum_{\gamma \in r^\vee \ol{Q}} e^{\lambda + n\gamma} q^{\abs{\lambda+n\gamma}^2/2n}.\tag{4.3}
\end{align*}
In order to study modular transformations of these theta forms, it is convenient to introduce new coordinates in $\h^*$:
\begin{align*}
    (\tau, z, t)_u = 2\pi i \pr{-\tau \frac{\L_0}{u'} + z + tu'\delta}, \quad \tau,t\in \C, z\in \ol{\h}^*.\tag{4.4}
\end{align*}
The coordinates, introduced in Section 1 are related to the new coordinates by
\begin{align*}
    (\tau,z,t) = \pr{u'\tau, z,\frac{t}{u'}}_u. \tag{4.5}
\end{align*}
In the new coordinates the theta form (4.3) looks as follows:
\begin{align*}
    \Theta_\lambda^{\#(u)}\pr{(\tau,z,t)_u} = e^{2\pi i n t} \sum_{\gamma \in r^\vee \ol{Q}} e^{2\pi i (\lambda + n\gamma \mid z)} q^{|\lambda+n\gamma|^2/2n}.\tag{4.6}
\end{align*}
Their modular transformations are given by the following formula (cf. \cite[Chapter 12]{K90}):
\begin{gather*}
    \Theta_{\lambda}^{\#(u)}\pr{\pr{-\frac{1}{\tau}, \frac{z}{\tau}, t - \frac{(z\mid z)}{2\tau}}_u} =  (-i\tau)^{\ell/2}\pr{\frac{n}{r^\vee}}^{-\ell/2} \abs{\ol{P}^\vee/r^\vee \ol{Q}}^{-1/2} \\
    \quad \times \sum_{\mu \in \frac{n}{u'} \L_0 + \frac{1}{r^\vee} \ol{P}^\vee \text{ mod } (nr^\vee \ol{Q} + \C \delta)} e^{-\frac{2\pi i}{n} \pr{\ol{\lambda}\mid\ol{\mu}}} \Theta_\mu^{\#(u)}\pr{\pr{\tau, z,t}_u}.\tag{4.7}
\end{gather*}

\section{Modular transformation of subprincipal admissible characters}

Recall that the span of the function $R_\g = A_\rho$ is $SL_2(\Z)$-invariant (see \cite{K90}, Lemma 13.8). Hence by (3.6), in order to study $SL_2(\Z)$-invariance of the span of characters of subprincipal admissible $\g$-modules $L(\L,\g)$, we need to study $SL_2(\Z)$-invariance of the span of functions $\set{A_{\L + \rho}}_{\L \in P_+^k}$ (see (2.10) - (2.12)). It is given by the following theorem, central to the paper.
\vspace{0.3cm}\\
\noindent \textbf{Theorem 5.1.} Let $k= \frac{p}{u} - h^\vee \ \pr{=\frac{p}{u'r^\vee}-h^\vee}$ be a subprincipal admissible level (see (2.1)), and let $\L \in P_{+,y}^k$, where $y = t_\beta \ol{y} \in \Phi_{(u)}$, be a subprincipal admissible weight (see (2.10), (2.11)), so that
\begin{align*}
    \L + \rho = y\pr{\frac{p}{u}\L_0 + \ol{v}}, \quad \text{where $\ol{v}\in \ol{P}$.}
\end{align*}
Then we have the following modular transformation formula for the numerators $A_{\L+\rho}$ (see (3.5)):
\begin{align*}
    A_{\L +\rho}\pr{-\frac{1}{\tau}, \frac{z}{\tau}, t - \frac{(z\mid z)}{2\tau}} = \sum_{y'=t_{\beta'}\ol{y}' \in \Phi_{(u)}} \ \sum_{\L' = y'(\frac{p}{u}\L_0 + \ol{v}')-\rho \in P_{+,y'}^k}a(\L,\L')A_{\L'+\rho}(\tau,z,t),
\end{align*}
where
\begin{gather*}
    a(\L,\L') = (-i\tau)^{\ell/2}(up)^{-\ell/2}\abs{\ol{P}^\vee/r^\vee Q}^{-1/2} \eps(\ol{y} \ol{y}')\\
    \qquad \qquad \qquad \times \sum_{\ol{w}\in\ol{W}}\eps(\ol{w}) e^{-\frac{2\pi i u}{p}\pr{\ol{v}\mid \ol{w}(\ol{v}')}}e^{-2\pi i\pr{\frac{p}{u}(\beta\mid \beta') + (\ol{v}\mid \beta') + (\ol{v}'\mid \beta)}}
\end{gather*}

In order to prove this theorem, we express the function $A_{\L + \rho}$, given by (3.5), in terms of theta forms (4.3), and then use the modular transformation formula (4.7) of these theta forms.

Introduce the following notation:
\begin{align*}
    \lambda_{(\ol{v},\beta,\ol{w})} := u't_\beta \ol{w}\pr{\frac{p}{u}\L_0 + \frac{\ol{v}}{r^\vee}}, \quad\text{where $\ol{v},\beta \in \ol{P}^\vee, \ol{w} \in \ol{W}$.}\tag{5.1}
\end{align*}
From (1.9) we have
\begin{align*}
    \lambda_{(\ol{v},\beta,\ol{w})} \equiv \frac{p}{r^\vee}\L_0 + \frac{u' \ol{w}(\ol{v}) + p \beta}{r^\vee} \text{ mod }\C\delta.\tag{5.2}
\end{align*}
Let
\begin{align*}
    \widetilde{\Omega}_{(u)} := \frac{p}{r^\vee}\L_0 + \frac{1}{r^\vee}\ol{P}^\vee; \qquad \Omega_{(k)} := \frac{p}{r^\vee}\L_0 + \frac{1}{r^\vee}\pr{\ol{P}^\vee/up\ol{Q}}.\tag{5.3}
\end{align*}
Let $C_{(u)}^\# = \set{\lambda \in \h^* \mid \ar{\lambda,\gamma}\geq 0\text{ for all $\gamma \in S_{(u)}^{\#\vee}$}}$ be the Weyl chamber for the Weyl group $W_{(u)}^\#$.

The following lemma is very important for the proof of Theorem 5.1.
\vspace{0.3cm}\\
\noindent\textbf{Lemma 5.1.} For $\mu\in \widetilde{\Omega}_{(k)}$ there exist $\ol{v}\in r^\vee \ol{P}$, $\beta\in r^\vee \ol{P}^\vee$ and $w\in \ol{W}$, such that the following two properties hold:
\begin{gather*}
    \mu \equiv \lambda_{(\ol{v},\beta,\ol{w}) }\text{ mod }\C\delta,\tag{5.4}\\
    \frac{p}{u'}\L_0 + \ol{v} \in C_{(u)}^\#.\tag{5.5}
\end{gather*}
\noindent{Proof.} Let
\begin{align*}
    I_\ell \text{ (resp. $I_s$)} = \set{1\leq i \leq \ell \mid \alpha_i \text{ is a long (resp. short) root}},
\end{align*}
and note that
\begin{align*}
    \ol{\L}_i^\vee = \ol{\L}_i \text{ (resp. $r^\vee \ol{\L}_i$) if $i\in I_\ell$ (resp. $I_s$)}.\tag{5.6}
\end{align*}
Since, by the conditions on $\mu$, $r^\vee \mu \in p\L_0 + \ol{P}^\vee$, we have
\begin{align*}
    r^\vee \mu = p \L_0 + \sum_{i=1}^\ell n_i\ol{\L}_i^\vee, \text{ where $n_i\in \Z$.}\tag{5.7}
\end{align*}
Since $\gcd(p,u) = 1 = \gcd(p,u')$, each $n_i$ decomposes as
\begin{align*}
    n_i = u'r^\vee n_i' \text{ (resp. $u'n_i'$)} + pn_i''\text{ if $i\in I_\ell$ (resp. $I_s$), where $n_i', n_i''\in\Z$.}\tag{5.8}
\end{align*}
Substituting (5.6) and (5.8) in (5.7), it becomes
\begin{align*}
    r^\vee \mu = u'v' + p \beta',\tag{5.9}
\end{align*}
where 
\begin{align*}
    v' = \frac{p}{u'}\L_0 + r^\vee\sum_{i=1}^\ell n_i'\ol{\L_i}, \quad \beta' = \sum_{i=1}^\ell n_i''\ol{\L}_i^\vee.\tag{5.10}
\end{align*}
We choose $w= t_{u'\xi}\ol{w} \in W_{(u)}^\#$, where $\xi \in r^\vee \ol{Q}$ and $\ol{w} \in \ol{W}$, such that
\begin{align*}
    v := w^{-1}v' \text{ lies in the Weyl chamber $C_{(u)}^\#$ for $W_{(u)}^\#$.}\tag{5.11}
\end{align*}
Note that
\begin{align*}
    v' = w(v) \equiv \ol{w}(v) + p \xi \text{ mod } \C\delta.\tag{5.12}
\end{align*}
Plugging (5.12) in (5.9), we obtain
\begin{align*}
    r^\vee \mu = u'v' + p\beta' \equiv u'\ol{w}(v) + p \beta \text{ mod }\C \delta, \text{ where $\beta = \beta' + u'\xi$}.\tag{5.13}
\end{align*}
Plugging $v = \frac{p}{u'}\L_0 + \ol{v}$ in (5.12), we obtain 
\begin{align*}
    \frac{p}{u'}\L_0 + \ol{v}' \equiv \frac{p}{u'}\L_0 + \ol{w}(\ol{v}) + p\xi \text{ mod }\C\delta,
\end{align*}
hence
\begin{align*}
    \ol{v} \equiv \ol{w}^{-1}(\ol{v}') - p \ol{w}' (\xi) \in r^\vee\ol{P},\tag{5.14}
\end{align*}
since $\ol{w}^{-1}\pr{\ol{v}'}\in r^\vee \ol{P}$ and $r^\vee \ol{Q} \subseteq r^\vee \ol{P}.$ By (5.13) and (5.9), we obtain
\begin{align*}
    r^\vee\mu \equiv u'\ol{w}\pr{\frac{p}{u'}\L_0 + \ol{v}}+p\beta \equiv u't_\beta \ol{w} \pr{\frac{p}{u'}\L_0 + \ol{v}} \text{ mod }\C\delta.
\end{align*}
Hence, mod $\C\delta$,
\begin{align*}
    \mu \equiv u't_\beta \ol{w} \pr{\frac{p}{u'}\L_0 + \frac{\ol{v}}{r^\vee}} = \lambda_{(\ol{v},\beta,\ol{w})}.\tag*{\qed}
\end{align*}
\newcommand{\td}{\widetilde}
\vspace{0.4cm}\\
\indent Introduce the following sets:
\begin{align*}
    \td{T}_{(k)}' \ (\text{resp. }T_{(k)}') = \big{\{}(\ol{v},\beta,\ol{w})\mid \ol{v}\in r^\vee \ol{P}\text{ is such that }\frac{p}{u'}\L_0 + \ol{v}\in C_{(u)}^\#, \\\beta\in\ol{P}^\vee \ (\text{resp. } \ol{P}^\vee/u\ol{Q}), \ol{w}\in\ol{W}\big{\}}.
\end{align*}
Denote by $\td{T}_{(k)}^\mathrm{reg}$ (resp. ${T}_{(k)}^\mathrm{reg}$) the subset of triples in $\td{T}_{(k)}'$ (resp. $T_{(k)}'$), such that $\ar{\ol{v},\alpha^\vee}\not\in pr^\vee\Z$ for all $\alpha^\vee\in\ol{\Delta}^\vee$, and let ${\td{T}^\mathrm{sing'}}_{(k)} = \td{T}_{(k)}' \setminus \td{T}_{(k)}^\mathrm{reg}$ and ${{T}^\mathrm{sing'}}_{(k)} = {T}_{(k)}' \setminus {T}_{(k)}^\mathrm{reg}$.
\vspace{0.4cm}\\
\textbf{Lemma 5.2.} $\td{\Omega}_{(k)} = \set{\lambda_{(\ol{v},\beta,\ol{w})} \mid(\ol{v},\beta,\ol{w}) \in \td{T}_{(k)}'}$, and similarly for $\Omega_{(k)}$.
\begin{proof}
    The inclusion $\subseteq $ follows from Lemma 5.1, while $\supseteq $ follows from
    \begin{align*}
        r^\vee \ol{P}\subseteq \ol{P}^\vee.\tag{5.15}
    \end{align*}
\end{proof}
Lemma 5.2 implies that the map
\begin{align*}
    (\ol{v},\beta,\ol{w})\> \lambda_{  (\ol{v},\beta,\ol{w})}\tag{5.16}
\end{align*}
gives bijections
\begin{align*}
    \td{T}_{(k)}^\mathrm{reg} \xrightarrow[]{\sim} \td{\Omega}_{(k)}^\mathrm{reg} \quad\text{and}\quad {T}_{(k)}^\mathrm{reg} \xrightarrow[]{\sim} {\Omega}_{(k)}^\mathrm{reg}.\tag{5.17}
\end{align*}
By Lemma 5.2, we can also choose subsets $\td{T}^\mathrm{sing}\subseteq \td{T}^\mathrm{sing'}$ and ${T}^\mathrm{sing}\subseteq {T}^\mathrm{sing'}$ such that the map (5.16) includes bijections 
\begin{align*}
    \td{\Omega}_{(k)}^\mathrm{sing} \xrightarrow[]{\sim} \td{T}_{(k)}^\mathrm{sing} \quad\text{and}\quad {\Omega}_{(k)}^\mathrm{sing} \xrightarrow[]{\sim} {T}_{(k)}^\mathrm{sing}.\tag{5.18},
\end{align*}
such that $(\ol{v},\beta,\ol{w})\in \td{T}_{(k)}^\mathrm{sing} $ for some $\ol{w}\in\ol{W}$ implies that $(\ol{v},\beta,\ol{w})\in \td{T}_{(k)}^\mathrm{sing} $ for all $\ol{w}\in\ol{W}$.

Also, it is straightforward to check the following,
\begin{align*}
    \pr{\ol{w}(\ol{v})\mid \beta} - \pr{\ol{v}\mid \beta}\in\Z \text{ if }\ol{w}\in\ol{W}, \ol{v}\in\ol{P},\beta\in\ol{P}^\vee,\tag{5.19}
\end{align*}
and
\begin{align*}
    \frac{r^\vee}{u'p}\pr{\ol{\lambda}_{\pr{r^\vee\ol{v},\beta,\ol{w}}} \mid \ol{\lambda}_{\pr{r^\vee\ol{v}',\beta',\ol{w}'}}} &\equiv \frac{u}{p}\pr{\ol{w}(\ol{v})\mid \ol{w}'(\ol{v}')} \\
    &\quad + \pr{\ol{v}\mid \beta'} + \pr{\ol{v}'\mid \beta} + \frac{p}{u}\pr{\beta\mid\beta'}\text{ mod }\Z
\end{align*}
for $\ol{v},\ol{v}'\in\ol{P}, \beta,\beta'\in \ol{P}^\vee, \ol{w},\ol{w}'\in\ol{W}$.
\vspace{0.4cm}\\
\noindent \textbf{Lemma 5.3.} Let $\ol{v}\in\ol{P}, (r^\vee\ol{v}',\beta',\ol{w}')\in T_{(k)}^\mathrm{sing}$. Then
\begin{align*}
    \sum_{w\in\ol{W}}\eps(\ol{w}) e^{-\frac{2\pi i u}{p}\pr{\ol{v}\mid \ol{w}(\ol{v}')}} = 0.\tag{5.21}
\end{align*}
\begin{proof}
    By definition of $T_{(k)}^\mathrm{sing}$, there exists $\alpha\in\ol{\Delta}$, such that $\pr{r^\vee \ol{v}'\mid \alpha^\vee} \in pr^\vee \Z$. But then, for any $\ol{w}\in\ol{W}$ we have:
    \begin{align*}
        \pr{\ol{w}^{-1}(\ol{v}) \mid r_\alpha(\ol{v}')}- \pr{w^{-1}(\ol{v})\mid \ol{v}'} = \pr{\ol{v}' \mid \alpha^\vee} \pr{\ol{w}^{-1}(\ol{v})\mid \alpha}\in\frac{p}{r^\vee}\Z,
    \end{align*}
    since $\pr{\ol{v}'\mid \alpha^\vee}\in p \Z$ and $\pr{\ol{w}^{-1}(\ol{v})\mid\alpha}\in\frac{1}{r^\vee}\Z$. Hence
    \begin{align*}
        \frac{u}{p}\pr{\pr{w^{-1}(\ol{v})\mid r_\alpha(\ol{v}')} - \pr{w^{-1}(\ol{v})\mid \ol{v}'} }\in u'\Z,\tag{5.22}
    \end{align*}
    and we have
    \begin{align*}
        \sum_{\ol{w}\in\ol{W}}\eps(\ol{w}) e^{-2\pi i \frac{u}{p}\pr{\ol{v}\mid \ol{w}(\ol{v}')}} &= \sum_{w\in \ol{W}} \eps(\ol{w} r_\alpha) e^{-2\pi i \frac{u}{p}\pr{\ol{v}\mid \ol{w}r_\alpha(\ol{v}')}}\\
        &= -\sum_{\ol{w}\in \ol{W}} \eps(\ol{w}) e^{-2\pi i \frac{u}{p}\pr{\ol{w}^{-1}(\ol{v})\mid r_\alpha(\ol{v}')}}.
    \end{align*}
    Using (5.22), this is equal to
    \begin{align*}
         -\sum_{\ol{w}\in \ol{W}} \eps(\ol{w}) e^{-2\pi i \frac{u}{p}\pr{\ol{w}^{-1}(\ol{v})\mid \ol{v}'}} = -\sum_{\ol{w}\in \ol{W}} \eps(\ol{w}) e^{-2\pi i \frac{u}{p}\pr{\ol{v}\mid \ol{w}(\ol{v}')}} ,
    \end{align*}
    proving (5.21).
\end{proof}
\textbf{Lemma 5.4.} The map
\begin{align*}
    \pr{r^\vee \ol{v},\beta,\ol{w}} \> \pr{\L = \frac{p}{u}\L_0 + \ol{v} - \rho, \beta, \ol{w}}
\end{align*}
gives an isomorphism
\begin{align*}
    T_{(k)}^\mathrm{reg} \xrightarrow[]{\sim} P_{+,e}^k \times \ol{P}^\vee/u\ol{Q} \times \ol{W}.
\end{align*}
\textit{Proof.} By definition of $T_{(k)}^\mathrm{reg} $, we have $\pr{r^\vee\ol{v},\beta,\ol{w}} \in T_{(k)}^\mathrm{reg} $ if and only if
\begin{enumerate}
    \item[(i)] $\ol{v}\in\ol{P}$ is such that $r^\vee \pr{\frac{p}{u} \L_0 + \ol{v}}$ lies in the open Weyl chamber for $W_{(u)}^\#$,
    \item[(ii)] $\beta \in \ol{P}^\vee/u\ol{Q}$,
    \item[(iii)] $\ol{w}\in\ol{W}$.
\end{enumerate}
But condition (i) is equivalent to
\begin{enumerate}
    \item[(i')] $\ol{v}\in \ol{P}$ is such that $\frac{p}{u}\L_0 + \ol{v}$ lies in the open Weyl chamber for $W_{(u)}^\#$, 
\end{enumerate}
which is equivalent to
\begin{enumerate}
    \item[(i'')] $\L := \frac{p}{u}\L_0 + \ol{v} - \rho \in P_{+,e}^k$,
\end{enumerate}
proving the lemma.\qed
\vspace{0.4cm}\\
\noindent\textbf{Proposition 5.1.} Let $\L \in P_{+,y}^k$ and $y\in t_\beta \ol{y}\in \Phi_{(u)}$ be such that
\begin{align*}
    \L + \rho = y\pr{\frac{p}{u}\L_0 + \ol{v}},\text{ where $\ol{v} \in \ol{P}$.}
\end{align*}
Then
\begin{align*}
    A_{\L+\rho}(\tau, z,t) = \eps(\ol{y})\sum_{w\in \ol{W}} \eps(\ol{w}) \Theta^{\#(u)}_{\lambda_{(r^\vee \ol{v}, \beta,\ol{w})}}\pr{\pr{\tau, \frac{z}{u'}, \frac{t}{u'^2}}_u}.\tag{5.23}
\end{align*}
\begin{proof}
    Recall the definition (3.5) of $A_{\L + \rho}$, and let there $w = y w' y^{-1}$, to obtain
    \begin{align*}
        A_{\L + \rho} = t_{\beta}\ol{y}\pr{q^{m_\L}\sum_{w\in W^\#_{(u)}} \eps(w) e^{w\pr{\frac{p}{u}\L_0 + \ol{v}}}}.
    \end{align*}
    Since $W_{(u)}^\# = \ol{W} \ltimes t_{u\ol{Q}}$ by (1.13), this formula can be rewritten as (since $\ol{y}t_\gamma=t_{\ol{y}(\gamma)}\ol{y}$):
    \begin{align*}
        A_{\L + \rho} = t_\beta \pr{q^{m_\L} \sum_{\ol{w}\in\ol{W}} \eps(\ol{w}) \sum_{\gamma \in u\ol{Q}} t_{\ol{y}(\gamma)}\ol{y} \pr{e^{\ol{w}\pr{\frac{p}{u}\L_0 + \ol{v}}}}}.
    \end{align*}
    Replacing $\ol{y}(\gamma)$ by $\gamma$ and $\ol{y}\ol{w}$ by $\ol{w}$, this can be rewritten as
    \begin{align*}
        A_{\L + \rho} = \eps(\ol{y})q^{m_\L} \sum_{\ol{w}\in\ol{W}} \eps(\ol{w}) \sum_{\gamma \in u\ol{Q}} t_{\beta}t_\gamma \pr{e^{\ol{w}\pr{\frac{p}{u}\L_0 + \ol{v}}}}.
    \end{align*}
    Writing this formula in coordinates $(\tau,z,t)$, using the new coordinates (4.4) and (4.5), we have
    \begin{align*}
        A_{\L + \rho}(\tau,z,t) = \eps(\ol{y}) q^{m_\L}\sum_{\ol{w}\in\ol{W}}\eps(\ol{w})\sum_{\gamma\in r^\vee\ol{Q}}t_{u'\gamma}\pr{e^{t_\beta \ol{y}\pr{\frac{p}{u}\L_0 + \ol{v}}}} \pr{u'\pr{\tau,\frac{z}{u'},\frac{t}{u'^2}}_u}.
    \end{align*}
    Using that $\lambda_{(r^\vee\ol{v},\beta,\ol{w})} = u't_\beta \ol{w}\pr{\frac{p}{u}\L_0 + \ol{v}}$ (see (5.11)), we obtain (5.23) using (4.3).
\end{proof}

\begin{proof}[Proof of Theorem 5.1.]
    Recall formula (4.7) for the modular transformation of the theta form $\Theta_\lambda^{\#(u)}$. Using (5.17), we deduce from it:
    \begin{align*}
        \Theta_{\lambda_{(r^\vee\ol{v},\beta,\ol{w})}}^{\#(u)}\pr{\pr{-\frac1\tau, \frac{z}{u'\tau}, \frac{t}{u'^2} - \frac{(z\mid z)}{2u'^2 \tau}}_u} = (-i\tau)^{\ell/2}(up)^{-\ell/2}\abs{\ol{P}^\vee/r^\vee\ol{Q}}^{-1/2}\\
        \times \sum_{(r^\vee\ol{v}',\beta',\ol{w}') \in T_{(k)}}e^{-\frac{2\pi i r^\vee}{u'p}\pr{\ol{\lambda}_{(r^\vee\ol{v},\beta,\ol{w})} \mid \ol{\lambda}_{(r^\vee\ol{v}',\beta',\ol{w}')}}}\Theta_{\lambda_{(r^\vee\ol{v}',\beta',\ol{w}')}}^{\# (u)}\pr{\pr{\tau, \frac{z}{u'}, \frac{t}{u'^2}}_u}.
    \end{align*}
    Using (5.19) and (5.20) we can rewrite this as follows:
    \begin{align*}
        \Theta_{\lambda_{(r^\vee\ol{v},\beta,\ol{w})}}^{\#(u)}\pr{\pr{-\frac1\tau, \frac{z}{u'\tau}, \frac{t}{u'^2} - \frac{(z\mid z)}{2u'^2 \tau}}_u} = (-i\tau)^{\ell/2}(up)^{-\ell/2}\abs{\ol{P}^\vee/r^\vee\ol{Q}}^{-1/2}\\
        \times \sum_{(r^\vee\ol{v}',\beta',\ol{w}') \in T_{(k)}}e^{-\frac{2\pi i u}{p}\pr{\ol{v}\mid \ol{w}^{-1}\ol{w}'(\ol{v}')}}e^{-2\pi i ((\ol{v}\mid \beta') + (\ol{v}'\mid \beta) + \frac{p}{u}(\beta\mid\beta'))}\\\times\Theta_{\lambda_{(r^\vee\ol{v}',\beta',\ol{w}')}}^{\# (u)}\pr{\pr{\tau, \frac{z}{u'}, \frac{t}{u'^2}}_u}.\tag{5.24}
    \end{align*}
    Using (5.24), we obtain from (5.23) the following for $\L = t_\beta \ol{y} \pr{\frac{p}{u}\L_0 + \ol{v}}-\rho \in P_{+,y}^k$, where $t_\beta \ol{y}\in \Phi_{(u)}$:
    \begin{align*}
        A_{\L + \rho}\pr{-\frac1\tau, \frac{z}{u'\tau}, \frac{t}{u'^2} - \frac{(z\mid z)}{2u'^2 \tau}} = (-i\tau)^{\ell/2}(up)^{-\ell/2}\abs{\ol{P}^\vee/r^\vee\ol{Q}}^{-1/2}\eps(\ol{y})\\
        \times \sum_{\ol{w}\in\ol{W}} \eps(\ol{w})\sum_{(r^\vee\ol{v}',\beta',\ol{w}') \in T_{(k)}}e^{-\frac{2\pi i u}{p}\pr{\ol{v}\mid \ol{w}^{-1}\ol{w}'(\ol{v}')}}e^{-2\pi i ((\ol{v}\mid \beta') + (\ol{v}'\mid \beta) + \frac{p}{u}(\beta\mid\beta'))}\\\times \Theta_{\lambda_{(r^\vee\ol{v}',\beta',\ol{w}')}}^{\# (u)}\pr{\pr{\tau, \frac{z}{u'}, \frac{t}{u'^2}}_u}.\tag{5.25}
    \end{align*}
    Letting $\ol{w} = \ol{w}'\ol{w}''^{-1}$ in the second line of (5.25), we can rewrite it as
    \begin{align*}
        \sum_{\ol{w}''\in\ol{W}} \sum_{(r^\vee\ol{v}',\beta',\ol{w}') \in T_{(k)}}\eps(\ol{w}') \eps(\ol{w}'')e^{-\frac{2\pi i u}{p}\pr{\ol{v}\mid \ol{w}''(\ol{v}')}}e^{-2\pi i ((\ol{v}\mid \beta') + (\ol{v}'\mid \beta) + \frac{p}{u}(\beta\mid\beta'))}\\\times \Theta_{\lambda_{(r^\vee\ol{v}',\beta',\ol{w}')}}^{\# (u)}\pr{\pr{\tau, \frac{z}{u'}, \frac{t}{u'^2}}_u}.\tag{5.26}
    \end{align*}
    The interior sum in the expression (5.26) can be rewritten as
    \begin{align*}
        \sum_{(r^\vee\ol{v}',\beta',\ol{w}') \in T_{(k)}} A_{v,v'} \eps(\ol{w}')e^{-2\pi i ((\ol{v}\mid \beta') + (\ol{v}'\mid \beta) + \frac{p}{u}(\beta\mid\beta'))}\Theta_{\lambda_{(r^\vee\ol{v}',\beta',\ol{w}')}}^{\# (u)}\pr{\pr{\tau, \frac{z}{u'}, \frac{t}{u'^2}}_u}.
    \end{align*}
    But $T_{(k)} =T_{(k)}^\mathrm{reg} \cup T_{(k)}^\mathrm{sing}$, and $A_{v,v'} = 0$ if $(r^\vee\ol{v}',\beta',\ol{w}')\in T_{(k)}^\mathrm{sing}$, by Lemma 5.3. Hence, using Lemma 5.4, we can rewrite (5.26) as
    \begin{align*}
        \sum_{(r^\vee\ol{v}',\beta',\ol{w}') \in T_{(k)}}\pr{\sum_{\ol{w}''\in\ol{W}} \eps(\ol{w}'')e^{-2\pi i \frac{u}{p} \pr{\ol{v}'\mid \ol{w}'' (\ol{v}')}} } \\
        \times \eps(\ol{w}')e^{-2\pi i ((\ol{v}\mid \beta') + (\ol{v}'\mid \beta) + \frac{p}{u}(\beta\mid\beta'))}\Theta_{\lambda_{(r^\vee\ol{v}',\beta',\ol{w}')}}^{\# (u)}\pr{\pr{\tau, \frac{z}{u'}, \frac{t}{u'^2}}_u}.\tag{5.27}
    \end{align*}
    Using Lemma 5.4, we can rewrite (5.27) as
    \begin{align*}
        \sum_{\ol{v}'\in\ol{P}}\sum_{\beta'\in\ol{P}^\vee/u\ol{Q}}\pr{\sum_{\ol{w}''\in\ol{W}} \eps(w'')e^{-2\pi i \frac{u}{p} \pr{\ol{v}'\mid \ol{w}'' (\ol{v}')}} } e^{-2\pi i ((\ol{v}\mid \beta') + (\ol{v}'\mid \beta) + \frac{p}{u}(\beta\mid\beta'))}\\
        \times \sum_{\ol{w}'\in\ol{W}}\eps(w') \Theta_{\lambda_{(r^\vee\ol{v}',\beta',\ol{w}')}}^{\# (u)}\pr{\pr{\tau, \frac{z}{u'}, \frac{t}{u'^2}}_u},\tag{5.28}
    \end{align*}
    where $\ol{v}'\in\ol{P}$ is such that $\frac{p}{u}\L_0 + \ol{v}' = \L_e' + \rho$, where $\L_e' \in P_{+,e}^k$. But, by Proposition 5.1, the second line in (5.28) is equal to $\eps(y') A_{\L' + \rho}(\tau,z,t)$, where $y' = t_{\beta'} \ol{y}'\in\Phi_{(u)}$, and
    \begin{align*}
        \L' + \rho := t_{\beta'}\ol{y}' \pr{\L_e' + \rho} \in P_{+,y'}^k.
    \end{align*}
    This completes the proof of Theorem 5.1.
\end{proof}

\section{Quantum Hamiltonian reduction $H_f(\L)$ for subprincipal admissible weights $\L$}

In \cite[Theorem 5]{KW25} we wrote down an explicit formula for the quantum Hamiltonian reduction $H_f(\L)$, associated to a non-zero nilpotent element $f\in \bar{\g}$ of a $\g$-module $L(\L)$ where $\L$ is a quasidominant principal admissible weight, which immediately implies conditions of non-vanishing of $H_f(\L)$, and the modular invariance of the vertex algebras $\td{W}_k(\bar{\g},f)$ for principal admissible levels $k$. In this section we do the same for quasidominant subprincipal admissible $\L$ of level $k$ and the vertex algebras $\td{W}_k(\bar{\g},f)$ for subprincipal admissible level $k$.
\vspace{0.4cm}\\
\textbf{Theorem 6.1.} Let $k= \frac{p}{u}-h^\vee$, where $u=u'r^\vee$ be a subprincipal admissible level. Let $\L \in P_{+,e}^k$ be a quasidominant subprincipal admissible weight of level $k$ for $\g$ (see Theorem 2.1), and let $(e,x,f)$ be an $\mathfrak{sl}_2$-triple in $\ol{\g}$. Then the normalized character of $H_f(\L)$ of the quantum Hamiltonian reduction of the $\g$-module $L(\L)$ is given by the following formula ($z\in \ol{\h}^f$):
\begin{align*}
    \mathrm{ch}\ H_f(\L) (\tau,z) = q^a e^{2\pi i \pr{\ol{\rho} - \ol{\rho}_{\g_0}\mid z}} \frac{\eta(u'\tau)^\ell \eta(u\tau)^{\frac{N-r^\vee\ell}{r^\vee-1}}}{\eta(\tau)^\ell}A(\tau,z)B(\tau,z)C(\tau,z)D(\tau,z),
\end{align*}
where
\begin{align*}
    A(\tau,z) &= \frac{\prod_{\alpha\in \ol\Delta_+,\alpha(x)>0}\pr{1-q^{\alpha(x)}e^{-2\pi i \alpha(z)}}}{\prod_{n=1}^\infty \pr{\prod_{\alpha\in \ol{\Delta}_+,\alpha(x)=1/2}\pr{1-q^{n-1/2}e^{2\pi i \alpha(z)}}\prod_{\alpha\in \ol{\Delta},\alpha(x)=0}\pr{1-q^{n}e^{2\pi i \alpha(z)}}}},\\
    B(\tau,z) &= \prod_{n=1}^\infty \pr{\prod_{\alpha\in \ol{\Delta}_\ell^+}\pr{1-q^{un-\alpha(x)}e^{2\pi i \alpha(z)}} \prod_{\alpha\in \ol{\Delta}_s^+} \pr{1-q^{u'n-\alpha(x)}e^{2\pi i \alpha(z)}}}\\
    C(\tau,z) &= \prod_{n=1}^\infty \pr{\prod_{\alpha\in \ol{\Delta}_\ell^+}\pr{1-q^{un+\alpha(x)}e^{-2\pi i \alpha(z)}} \prod_{\alpha\in \ol{\Delta}_s^+} \pr{1-q^{u'n+\alpha(x)}e^{-2\pi i \alpha(z)}}}\\
    D(\tau,z) &= \mathrm{ch}\ L\pr{\L^{\# 0}, \g^{\#}_{(u)}}\pr{u'\tau, -\tau x + z, \frac{\tau}{2u'}\pr{x\mid x}},\\
    a &= \frac{1}{2hu'}\abs{u'\ol{\rho} - hx}^2 - \frac{
    u'\ell + u\frac{N-r^\vee \ell}{r^\vee -1} - \ell + \dim \ol{\g}_0 - \frac12 \dim \ol{\g}_{1/2}
    }{24}.
\end{align*}
\begin{proof}
    We compute the normalized Euler-Poincar\'{e} character of $H_f(\L)$ by the formula in \cite[Remark 3.1]{KRW03}. By \cite{A15}, the Poincar\'{e} character coincides with the character. Hence we have 
    \begin{align*}
        \mathrm{ch}H_f(\L) (\tau,z) = \frac{A_{\L+\rho}\pr{\tau, - \tau x + z, \frac{\tau}{2}\pr{x\mid x}}}{R^w(\tau,z)},\tag{6.1}
    \end{align*}
    where $R^w(\tau,z)$ is the normalized Weyl denominator for the universal $W$-algebra $W^k(\bar{\g},f)$ which is given by
    \begin{align*}
        R^w(\tau,z) = q^{\frac{1}{48}\pr{2\dim \ol{\g}_0 - \dim \ol{\g}_{1/2}}}e^{2\pi i \rho_{\g_0}(z)}Q^w\pr{\tau,z},\tag{6.2}
    \end{align*}
    where
    \begin{align*}
        Q^w(\tau,z) = \phi(q)^\ell \prod_{n=1}^\infty\Bigg{(}\prod_{\alpha\in\ol{\Delta}_+,\alpha(x)=0} \pr{1-e^{-(n-1)\delta-\alpha}}\pr{1-e^{-n\delta+\alpha}} \\\times \prod_{\alpha\in\ol{\Delta}_+,\alpha(x)=1/2} \pr{1-e^{-n\delta+\alpha}}\Bigg{)}(\tau,-\tau x+ z, 0),
    \end{align*}
    and $\phi(q) = \prod_{n=1}^\infty(1-q^n) = q^{1/24} \eta(\tau)$. Hence
    \begin{align*}
        Q^w(\tau,z) = \phi(q)^\ell \prod_{\alpha\in\ol{\Delta}_+,\alpha(x)=0} \pr{\pr{1-e^{-2\pi i \alpha(z)}} \prod_{n=1}^\infty \pr{1-q^ne^{-2\pi i \alpha(z)}} \pr{1-q^ne^{2\pi i \alpha(z)}} }\\
        \times \prod_{\alpha\in\ol{\Delta}_+,\alpha(x)=1/2}\prod_{n=1}^\infty \pr{1-q^{n-1/2}e^{2\pi i \alpha(z)}}\qquad \qquad
    \end{align*}
    Let $\g_{(u)}^\#$ be the affine subalgebra of the affine Lie algebra $\g^\#$ with the set of simple coroots $S_{(u)}^{\#\vee}$ (see (2.3)), and let $R_{(u)}^\#$ be its normalized Weyl denominator. It is given by
    \begin{align*}
        R_{(u)}^\# = q^{\frac{u'}{2h}\abs{\ol{\rho}}^2}Q_{(u)}^\#,\tag{6.4}
    \end{align*}
    where
    \begin{align*}
        Q_{(u)}^\# = e^{\rho_{(u)}^\#}\phi(q^{u'})^\ell \phi(q^u)^{\frac{N-r^\vee \ell}{r^\vee -1 }}\prod_{n=1}^\infty\Bigg{(}\prod_{\alpha\in\ol{\Delta}_\ell^+} \pr{1-q^{un}e^\alpha}\pr{1-q^{u(n-1)}e^{-\alpha}} \\\times\prod_{\alpha\in\ol{\Delta}_s^+} \pr{1-q^{u'n}e^\alpha}\pr{1-q^{u'(n-1)}e^{-\alpha}}\Bigg{)}.\tag{6.5}
    \end{align*}
    Since
    \begin{align*}
        \rho_{(u)}^\#\pr{\tau,-\tau x + z, \frac{\tau}{2}\pr{x\mid x}} = 2\pi i \tau \pr{\frac{h}{2u'}\pr{x\mid x} - \ol\rho(x)} + 2\pi i \rho(z),
    \end{align*}
    we obtain from (6.5):
    \begin{align*}
        Q_{(u)}^\#\pr{\tau,-\tau x + z, \frac{\tau}{2}\pr{x\mid x}}= q^{\frac{h}{2u'}\pr{x\mid x} - \ol\rho(x)}e^{2\pi i \ol{\rho}(z)}\phi\pr{q^{u'}}^\ell\\
        \times \phi(q^u)^{\frac{N-r^\vee \ell}{r^\vee -1 }}\prod_{n=1}^\infty\Bigg{(}\prod_{\alpha\in\ol{\Delta}_\ell^+} \pr{1-q^{un-\alpha(x)}e^{2\pi i \alpha(z)}}\\
        \times \pr{1-q^{u(n-1)+\alpha(x)}e^{-2\pi i \alpha(z)}}\prod_{\alpha\in\ol{\Delta}_s^+} \pr{1-q^{u'n-\alpha(x)}e^{2\pi i\alpha(z)}}\\
        \times \pr{1-q^{u'(n-1)+\alpha(x)}e^{-2\pi i\alpha(z)}}\Bigg{)}.\tag{6.6}        
    \end{align*}
    We can write (6.1) as a product of two factors:
    \begin{align*}
        \mathrm{I} = \frac{A_{\L+\rho}\pr{\tau, - \tau x + z, \frac{\tau}{2}\pr{x\mid x}}}{R^{\#}_{(u)}(\tau,-\tau x + z, \frac{\tau}{2}\pr{x\mid x})}\quad \text{and}\quad \mathrm{II} =  \frac{R^{\#}_{(u)}(\tau,-\tau x + z, \frac{\tau}{2}\pr{x\mid x})}{R^w(\tau,z)}.\tag{6.7}
    \end{align*}
    By Lemma 3.2(b) and (4.5), we have
    \begin{align*}
        \mathrm{I} = \mathrm{ch}\ L\pr{\L^{\# 0}, \g_{(u)}^\#}\pr{u'\tau, - \tau x + z, \frac{\tau}{2u'}\pr{x\mid x}}_u.\tag{6.8}
    \end{align*}
    Hence, by (6.1), (6.7) and (6.8), we have
    \begin{align*}
         \mathrm{ch} \ H_f(\L)(\tau,z)=\frac{R^{\#}_{(u)}(\tau,-\tau x + z, \frac{\tau}{2}\pr{x\mid x})}{R^w(\tau,z)}\ \mathrm{ch}\ L\pr{\L^{\# 0}, \g_{(u)}^\#}\pr{u'\tau, - \tau x + z, \frac{\tau}{2u'}\pr{x\mid x}}_u.\tag{6.9}
    \end{align*}
    Inserting (6.4), (6.5), (6.2), and (6.3) in this formula, and using (6.6), we complete the proof of the theorem.
\end{proof}
\noindent \textbf{Corollary 6.1} Let $k$ and $\L \in P_{+,e}^k$ be as in Theorem 6.1, and let $f$ be a non-zero nilpotent element of $\ol{\g}$ with Dynkin characteristic $2x$. Then
\begin{enumerate}
    \item[(a)] $H_f(\L)$ is zero if and only if there exist $n\in \Z_{\geq 1}$ and $\alpha \in \ol{\Delta}_+$, such that $\alpha(x) = nu$ for $\alpha \in \ol{\Delta}_\ell^+$, and $=nu'$ for $\alpha \in \ol{\Delta}_s^+$; equivalently, if and only if $u' \leq \theta_s(x)$.
    \item[(b)] If $f=f_\mathrm{min}$, then $H_f(\L)= W_k^\mathrm{min}(\ol{\g})$ is always non-zero.
    \item[(c)] If $f = f_\mathrm{pr}$, then $H_f(\L)$ is non-zero if and only if $u' \geq h^\#$, the Coxeter number of the affine Lie algebra $\g^\#$.
    \item[(d)] If $u'> \theta_s(x)$, then $H_f(k\L_0)$ is a modular invariant vertex algebra.
\end{enumerate}
\begin{proof}
    The proof of (a) is the same as of Corollary 1(a) from \cite{KW25}, using Lemma 5.15(ii) from \cite{A15}. Claims (b) and (c) follow from (a). Claim (d) follows from the fact that the function $D(\tau,0)$ is a modular function \cite[Chapter 13]{K90}.

    Note that Corollary 6.1 is equivalent to Theorem 5.16(ii) in \cite{A15}.
\end{proof}

\subsection*{Funding and/or Conflicts of Interests, Competing Interests} The authors have no conflict of interest to declare that are relevant to the content of the article.

\subsection*{Data Availability Statement} Data sharing is not applicable to this article as no datasets were generated or analyzed during the current study.



\begin{thebibliography}{9} 
\addtolength{\leftmargin}{0in} 
\setlength{\itemindent}{-0.5in}
\bibitem[Ara15]{A15} T. Arakawa. \emph{Associated varieties of modules over Kac-Moody algebras and $C_2$-cofiniteness of W-algebras}. IMRN(22), 2015,
pp. 11605–11666.
\bibitem[Ara16]{A16} T. Arakawa. \emph{Rationality of admissible affine vertex algebras in the category $\mathcal{O}$}. no. 1. Duke Math. J. 165, 2016, pp. 67–93.
\bibitem[Kac90]{K90} V.G. Kac. \emph{Infinite-dimensional Lie algebras}. Third edition. Cambridge University Press, 1990.
\bibitem[KRW03]{KRW03} V.G. Kac, S-S. Roan, and M. Wakimoto. \emph{Quantum reduction for affine superalgebras}. Commun. Math. Phys. 241, 2003, pp. 307-342
\bibitem[KW88]{KW88} V.G. Kac and M. Wakimoto. \emph{Modular invariant representations of infinite-dimensional Lie algebras and superalgebras}. Proc. Nat. Acad. Sci. USA 85, 1988, pp. 4956-4960.
\bibitem[KW89]{KW89} V.G. Kac and M. Wakimoto. \emph{Classification of modular invariant representations of affine algebras}. Vol. 7. Infinite-dimensional Lie algebras and groups (Luminy-Marseille, 1988), Advanced ser. in Math. Phys. World Scientific, 1989, pp. 138–177.
\bibitem[KW08]{KW08} V.G. Kac and M. Wakimoto. \emph{On rationality of W-algebras}. Vol. 7. Transform. Groups 13 no. 3-4. 2008, pp. 671–713.
\bibitem[KW25]{KW25} V.G. Kac and M. Wakimoto. \emph{On modular invariance of vertex operator algebras}. no. 2, Paper No. 44. Commun. Math. Phys. 406, 2025.
\end{thebibliography}
\end{document}